%% file: ms.tex
\documentclass[letterpaper, 10 pt, conference]{ieeeconf}  

\IEEEoverridecommandlockouts                              

\usepackage{graphics} 
\usepackage{epsfig} 
\usepackage{mathptmx} 
\usepackage{times} 
\usepackage{amsmath} 
\usepackage{amssymb}  

\usepackage{type1cm}  

\usepackage{bm}
\usepackage{xspace}
\usepackage{mathrsfs}
\usepackage{lettrine}
\usepackage[labelsep=period,font=bf]{caption}
\usepackage{exscale}

\usepackage{comment}

\usepackage[xspace]{ellipsis}

\usepackage{resizegather}

\usepackage{epstopdf}

\input{./LatexCommon/command}

\input{MacrosAndDefinitions}

\newcommand*{\myDots}{\ifmmode\mathellipsis\else.\kern-0.23em.\kern-0.23em.\fi} 

\newtheorem{proposition}{Proposition}
\newtheorem{theorem}{Theorem}
\newtheorem{corollary}{Corollary}

\includecomment{ExpandedCalculation}

\graphicspath{ {./figures/} }

\title{\LARGE \bf
Finite-Horizon Covariance Control of Linear Time-Varying Systems
}

\author{Maxim Goldshtein$^{1}$ and Panagiotis Tsiotras$^{2}$
\thanks{$^{1}$PhD student, Daniel Guggenheim School of Aerospace Engineering, 270 Ferst Dr., Atlanta, GA, 30332-0150, USA.
        Email:{maxg@gatech.edu}}%
\thanks{$^{2}$Professor, Daniel Guggenheim School of Aerospace Engineering, 270 Ferst Dr., Atlanta, GA, 30332-0150, USA.
        Email:{tsiotras@gatech.edu}}%
}

\begin{document}

\maketitle
\thispagestyle{empty}
\pagestyle{empty}

\begin{abstract}

We consider the problem of finite-horizon optimal control of a discrete linear time-varying system subject to a stochastic disturbance and fully observable state. The initial state of the system is drawn from a known Gaussian distribution, and the final state distribution is required to reach a given target Gaussian distribution, while minimizing the expected value of the control effort.
We derive the linear optimal control policy by first presenting an efficient solution for the diffusion-less case, and we then solve the case with diffusion by reformulating the system as a superposition of diffusion-less systems. This reformulation leads to a simple condition for the solution, which can be effectively solved using numerical methods.
We show that the resulting solution coincides with a LQG problem with particular terminal cost weight matrix. 
This fact provides an additional justification for using a linear in state controller. In addition, it allows an efficient iterative implementation of the controller.

\end{abstract}

\section{Introduction}

The work in this paper is aimed at solving the problem of the optimal steering of a discrete stochastic linear system, with a fully observable state, a known Gaussian distribution of the initial state, and a state and input-independent white-noise Gaussian diffusion.
The control task is to find the optimal input to steer the state of the system to a pre-specified target Gaussian distribution in a given time, while minimizing the expected value of the input signal $\ell_2$-norm.
Since the Gaussian distribution can be fully defined by its first two moments, this problem can be described as a finite-time optimal mean and covariance steering of a stochastic time-varying discrete linear system.
Unlike the classical LQG case~\cite{Astrom1970}, where the final state covariance appears as a by-product of the solution, here we are required to reach exactly the target covariance at the given final time.

The covariance steering problem is relevant to a wide range of control and path-planning applications, such as decentralized control of swarm robots~\cite{Shahrokhi2016}, closed-loop cooling~\cite{Vinante2008}, and other areas, where it is more natural to specify a distribution over the state rather than a fixed set of values.


The steady-state covariance control problem, also known as the Covariance Assignment problem, has been extensively studied for both continuous and discrete-time stochastic linear systems~\cite{Hotz1987,GrigoriadisSkelton1997,Xu1992,Collins1987}.
A finite-time optimal solution for the continuous case has been recently derived in~\cite{Chen.etal2016,Chen.etal2016a}, and~\cite{Chen2016Thesis}, with a connection to the problems of Shr{\"{o}}dinger bridges~\cite{Schrodinger1931} and the Optimal Mass Transfer~\cite{Kantorovich1942}.
In these works the authors showed that, if the diffusion term affects the system through all control input channels, the target probability can always be achieved in finite time, and the solution is given in state-feedback form.
A more general case, in which the control input and the diffusion channels are different, can be solved using a soft constraint on the target distribution (such as using the Wasserstein distance~\cite{Halder2016ACC}), or by numerical optimization methods~\cite{Bakolas2016ACC}.

The discrete finite-time case has been addressed in~\cite{Bakolas2016CDC}, in which the author used a relaxed formulation for the target covariance in order to facilitate its numerical solution. That is, the final states are expected to have a more concentrated distribution that does not necessarily follow the target Gaussian distribution.
In this paper we treat a similar problem as in~\cite{Bakolas2016CDC}, but we impose a hard equality constraint in the final distribution instead, so the relaxation imposed in~\cite{Bakolas2016CDC} is not needed.
In addition, the solution in~\cite{Bakolas2016CDC} is based on a non-linear convex programming with a large number of variables (approximately proportional to state size $\times$ inputs size $\times$ number of steps). The proposed method, on the other hand, requires only $n^2/2 + n$ decision variables (where $n$ is the state size), thus greatly reducing the required computations.

Another special case of linear discrete finite-time Gaussian stochastic systems was mentioned in~\cite{Beghi1996}, in which the author shows a relation between the relative entropy and the minimum energy LQG optimal control problems. The system discussed in~\cite{Beghi1996} has a full control authority and the disturbance matrix is invertible.
This paper extends the results presented in~\cite{Beghi1996} to a general linear system, In addition, the conditions for the solvability presented in~\cite{Beghi1996} follow naturally from the analysis presented here.

\subsubsection*{Main Contribution}

In this paper we first derive the minimum-control-effort optimal steering solution for linear time-varying discrete stochastic systems, subject to boundary conditions in terms of their Gaussian distribution.
The problem considered herein can be viewed as a subset of the problems presented in~\cite{Bakolas2016CDC}, but with a different solution formulation.
The paper provides necessary conditions for the existence of the solution, and proposes an efficient numerical scheme for attaining it.
The proposed solution assumes full-state observation, and provides an optimal policy that depends linearly on the state.
Such a linear dependence provides a natural retention\footnote{It is possible however that some other non-linear law will yield a better result, while retaining the Gaussian property of the state, but such a law, in general, will be very difficult to describe~\cite{Hamedani2001}.} for the Gaussian nature of the state, and has been proven to be the optimal policy in the case of linear continuous-time invariant systems with Gaussian state distribution~\cite{Brockett2007}. Furthermore, we show that the resulting controller coincides with solving a LQG~\cite{Astrom1970} problem, with the particular choice of the terminal cost weight matrix. 

The notation used throughout this paper is quite standard.
A unit matrix is denoted as $\I$, with a subscript indicating its dimensions, where relevant.
The notation $\Expectation{\cdot}$ denotes the expectation operator.
A random variable $x$ with normal distribution is denoted as $x \sim \normal{\mu, \Sigma}$, where $\mu$ is its mean, and $\Sigma$ its covariance matrix.
The trace of a square matrix is denoted by $\Trace{\cdot}$.
The positive-definiteness of the square matrix $R$ is denoted as $\posdef{R}$, and semi-definiteness is denoted as $\possemidef{R}$. A zero matrix with dimensions $m \times n$ is denoted as $0_{m \times n}$.
An $n \times n$  diagonal matrix with $(a_{1}, a_2, \ldots, a_n)$ on the diagonal is denoted as $\mathrm{diag}[a_1, a_2, \ldots, a_n]$.

The rest of this paper is organized as follows. In Section~\ref{secProblemStatement} we formulate the covariance control problem for discrete-time linear systems as a constrained minimization problem. In Section~\ref{secOptimalCovSteering} we provide an analytical solution to this problem.
In Section~\ref{secNumericalExample} we provide numerical examples to demonstrate the performance of the proposed method.
We also show that  the proposed solution recovers the same controller of the target covariance resulting from solving a LQG problem.
We finally conclude with a summary of the results and some suggestions for future work.

\section{Problem Statement}\label{secProblemStatement}

\subsection{Problem Formulation}\label{sec:secProblemStatement:ProblemFormulation}

Consider the discrete stochastic linear time-varying system
\begin{equation}\label{eqTheSystemDiscrete}
\x_{k+1} = \A_k\x_{k} + \B_{k}\u_{k} + \G_{k}\w_{k}. \quad k=0,1,...,N,
\end{equation}
where $\x \in \Reals{n}$ is the state, $\u \in \Reals{m}$ is the control input, and $\w \in \Reals{r}$  is a zero-mean white Gaussian noise with unit covariance, that is,
\begin{equation}\label{eqNoiseMdl1}
	\Expectation{\w_{k_1}\w_{k_2}^{\TR}}= \begin{cases}
	\I_r & k_1 = k_2, \\
	0  &  k_1 \neq k_2,
	\end{cases}
\end{equation}
and
\begin{equation}\label{eqNoiseMdl2}
	\Expectation{\x_{k_1}\w_{k_2}^{\TR}} = 0,  \qquad  0 \leq \k_1 \leq \k_2 \leq N.
\end{equation}

Our objective is to steer the trajectories of system~\eqref{eqTheSystemDiscrete} from a given initial Gaussian distribution having
mean $\Expectation{\x_{0}} = \muS$ and covariance $\covS$ to a final Gaussian distribution having mean $\Expectation{\x_{N+1}} = \muF$ and covariance $\covF$.
That is, we wish the initial and final states to be distributed according to
\begin{equation}\label{eqBoundaryConditions}
\x_0 \sim \normal{\muS,\covS}, \qquad \x_{N+1} \sim \normal{\muF,\covF},
\end{equation}
with $\muS,\covS,\muF,\covF$ given,
while minimizing the cost function
\begin{equation}\label{eqDiscreteCostFunction}
J(\u_0, ..., \u_N) = \Expectation{\Sum{k=0}{N}{ \u_{k}^\TR  \u_{k}}}.
\end{equation}

\subsection{System Dynamics}

For each step $k$, the system state $x_k$ can be explicitly calculated as follows.
Let $\A_{k_1, k_0}$, $\B_{k_1, k_0}$, $\G_{k_1, k_0}$ denote the transition matrices of the state, the input, and the diffusion term
from step $k_0$ to step $k_1+1$ ($k_1 > k_0$) as follows
\begin{subequations}
	\begin{align}\label{eqDiscreteTransitionMatricesDefinitions}
	\A_{k_1, k_0} &= \A_{k_1} \A_{k_1-1} \cdots \A_{k_0} ,\\
	\B_{k_1, k_0} &= \A_{k_1, k_0+1} \B_{k_0} ,\\
	\G_{k_1, k_0} &= \A_{k_1, k_0+1} \G_{k_0} ,
	\end{align}
\end{subequations}
and let, for simplicity,  $\A_{k, k} = \A_{k}, \B_{k, k} = \B_{k}, \G_{k, k} = \G_{k}$.
Let also $\U_{k_1, k_2}$ and $\W_{k_1, k_2}$ ($k_1 \leq k_2$) be the vectors
\begin{align}\label{eqDiscreteInputVectorsDefinitions}
\U_{k_1, k_2} = \bmat{\u_{k_1} \\ \u_{k_1+1} \\ \vdots \\  \u_{k_2}}, \qquad \W_{k_1, k_2} = \bmat{\w_{k_1} \\ \w_{k_1+1} \\ \vdots \\  \w_{k_2}},
\end{align}
and, for simplicity, let $\U_k \dfn \U_{0, k}$ and $\W_k \dfn \W_{0, k}$.
For convenience, define the matrices
\begin{subequations}
	\begin{align}\label{eqDiscreteTransitionMatricesDefinitionsCombined}
	\Bbar_{k_1, k_0} &\dfn \bmat{\B_{k_1, k_0} &  \B_{k_1, k_0+1}  & \cdots & \B_{k_1, k_1}}, \\
	\Gbar_{k_1, k_0} &\dfn \bmat{\G_{k_1, k_0} &  \G_{k_1, k_0+1}  & \cdots & \G_{k_1, k_1}},
	\end{align}
\end{subequations}
and let
\begin{equation}\label{eqTransMatricesShortcut}
	\Abar_k \dfn \A_{\k, 0}, \quad \Bbar_k \dfn \Bbar_{k, 0}, \quad \Gbar_k \dfn \Gbar_{k, 0}.
\end{equation}
The system state at step $k+1$ is given by
\begin{align}\label{eqDiscreteSystemEvaluationCalc}
\x_{k+1} &= \Abar_k\x_0 + \Bbar_k\U_k + \Gbar_k\W_k.
\end{align}

Since $\Expectation{\W_k} = 0$, the mean of the state obeys the expression
\begin{equation}\label{eqDiscreteStateStatisticsMean}
\mean_{k+1} \dfn \Expectation{\x_{k+1}}  = \Abar_k \mean_0 + \Bbar_k \Expectation{\U_k}.
\end{equation}
Defining now
\begin{equation}\label{eqUandxTildaDef}
\Utilde_k \dfn \U_k - \Expectation{\U_k}, \qquad   \xtilde_k \dfn \x_k - \mean_k,
\end{equation}
it follows that
\begin{equation}\label{eqXtildeDynamics}
\xtilde_{k+1} = \Abar_k\xtilde_0 + \Bbar_k\Utilde_k + \Gbar_k\W_k.
\end{equation}
The state covariance is given by\footnote{A causal state-feedback controller at step $k$ is independent of the diffusion term at step $k'$, with $k'\geq k$.}
\begin{align} \label{eqDiscreteStateStatisticsCov}
\mSigma_{k+1} &\dfn \Expectation{\xtilde_{k+1}\xtilde_{k+1}^{\TR}}   \notag\\
&=\Abar_k \mSigma_{0}\Abar_k^{\TR} +  \Bbar_k \Expectation{\Utilde_k \Utilde_k^{\TR}} \Bbar_k^{\TR} + \Gbar_k \Expectation{\W_k\W_k^{\TR}}\Gbar_k^{\TR} \notag\\
&\quad + \Bbar_k \Expectation{\Utilde_k \xtilde_0^{\TR}} \Abar_k^{\TR} + \Abar_k \Expectation{\xtilde_0\Utilde_k^{\TR}}\Bbar_k^{\TR} \notag\\
&\quad + \Bbar_k \Expectation{\Utilde_k \W_{k-1}^{\TR}} \Gbar_{k-1}^{\TR} + \Gbar_{k-1} \Expectation{\W_{k-1}\Utilde_k^{\TR}}\Bbar_k^{\TR}.
\end{align}
From~\eqref{eqUandxTildaDef}, the cost function~\eqref{eqDiscreteCostFunction} can be written as

\begin{align}\label{eqDiscreteCostFunctionVectirized}
J(\U_N) &= \Expectation{\U_N^{\TR}\U_N} = \underbrace{\Expectation{\U_N}^{\TR} \Expectation{\U_N}}_{\displaystyle{J_{\mu}}} + \underbrace{\Trace{\Expectation{\Utilde_N \Utilde_N^{\TR}}}}_{\displaystyle{J_{\Sigma}}}.
\end{align}

It will be assumed in this paper that the system~\eqref{eqTheSystemDiscrete} is \textit{controllable}, that is, if $G_k \equiv 0$, the reachable set at $k = {N+1}$ is $\Reals{n}$, that is, for any $\x_S\in\Reals{n}$ and $\x_F\in\Reals{n}$, there exist a set of controls $\brf{u_k}_{k=0}^N$ that brings the state from $x_0 = x_S$ to $x_{N+1} = x_F$.
From~\eqref{eqDiscreteSystemEvaluationCalc} it is straightforward to conclude that system~\eqref{eqTheSystemDiscrete} is controllable if and only if $\Bbar_N$ is full row rank.

\section{Optimal Covariance Steering}\label{secOptimalCovSteering}

As seen from~\eqref{eqDiscreteStateStatisticsMean},~\eqref{eqDiscreteStateStatisticsCov} and~\eqref{eqDiscreteCostFunctionVectirized}, the problem of steering the mean and the covariance can be separated into two independent subproblems: finding an optimal $\Expectation{\U_N}$ that minimizes $J_{\mu}$ satisfying the mean constraint~\eqref{eqDiscreteStateStatisticsMean} and the boundary condition~\eqref{eqBoundaryConditions}, and finding an optimal $\Utilde_N$ that minimizes $J_{\Sigma}$ satisfying the covariance constraint~\eqref{eqDiscreteStateStatisticsCov} and the boundary condition~\eqref{eqBoundaryConditions}. This section presents an analytical solution to both problems.

\subsection{Steering the Mean}\label{sec:OptimalSteering:Mean}

Since the dynamics of the state mean are governed by~\eqref{eqDiscreteStateStatisticsMean}, and the cost function that is influenced by the mean is given in~\eqref{eqDiscreteCostFunctionVectirized}, the optimal solution for $\Expectation{\U_N}$ will not influence the covariance part of the solution. The solution for the mean steering is well known in the literature, and is given below for the sake of completeness.

\vspace{1em}
\begin{proposition}\label{lemma-steering-the-mean}
	Given the controllable system~\eqref{eqTheSystemDiscrete}, the optimal control $\Expectation{\optU_N}$ that minimizes the cost
	\begin{equation*}
		J_{\mu} = \Expectation{\U_N}^{\TR} \Expectation{\U_N} = \sum_{k=0}^N \Expectation{u_k}^{\TR}  \Expectation{u_k},
	\end{equation*}	
	subject to the constraint
	\begin{equation}  \label{eq786}
	\Abar_N \muS + \Bbar_N \Expectation{\U_N} = \muF,
	\end{equation}
	is given by
	\begin{equation}\label{eqDiscreteOptimalSteeringMean}
	\Expectation{\optU_N} = \Bbar_N^{\TR}\invb{\Bbar_N \Bbar_N^{\TR}}\br{\muF - \Abar_N\muS}.
	\end{equation}
\end{proposition}
\vspace{1em}

\begin{proof}
	A reformulation of the constraint~\eqref{eq786} yields the algebraic linear system
	\begin{equation*}
		 \Bbar_N \Expectation{\U_N} = \muF - \Abar_N \muS.
	\end{equation*}
	Since the system~\eqref{eqTheSystemDiscrete} is controllable, $\Bbar_N$ has full row rank, and the solution to the above linear system, that minimizes the quadratic norm of $\Expectation{\U_N}$, is given by~\eqref{eqDiscreteOptimalSteeringMean}~\cite{Corless2003}.
\end{proof}
\vspace{1em}

Now that we have the mean steering solution, the rest of the paper will concentrate on solving the covariance steering problem, using the deviation-from-mean dynamics given by~\eqref{eqXtildeDynamics}, and the covariance-part cost $\displaystyle{J_{\Sigma}}$ given in~\eqref{eqDiscreteCostFunctionVectirized}. 
For simplicity, we will assume that the original system has zero-mean constraints for the initial and final states.

\subsection{Steering the Covariance}

In this section we present the covariance steering controller by first deriving a necessary condition for the solution, and then presenting a numerical scheme to find a controller that satisfies these necessary conditions.

To this end, assume a controller of the form
\begin{equation}\label{eqDiscreteLinearPolicy}
\Utilde_N = \L \xtilde_0,
\end{equation}
where $\L \in \Reals{(N r) \times n}$.
The covariance-related part of the cost function~\eqref{eqDiscreteCostFunctionVectirized} can now be rewritten as:
\begin{equation}\label{eqCovarianceCostFunctionLinearController}
J_{\Sigma} = \Trace{\Expectation{\Utilde_N \Utilde_N^\TR}} = \Trace{\L \covS \L^{\TR}},
\end{equation}  	
where $\covS = \Expectation{\xtilde_{0}\xtilde_{0}^{\TR}}$.

\vspace{1em}
\subsubsection{Diffusion-less Case}

Suppose that $\G_k = 0$ for all $k \in [0, N]$ in \eqref{eqTheSystemDiscrete}.
In this case, the final state covariance~\eqref{eqDiscreteStateStatisticsCov} becomes
\begin{align}\label{eqDiscreteStateStatisticsCovNoDiffusionAndDemeaned}
\covF = \mSigma_{N+1} &= \Abar_N \mSigma_{0}\Abar_N^{\TR} +  \Bbar_N \Expectation{\Utilde_N \Utilde_N^{\TR}} \Bbar_N^{\TR}  \notag\\
&\quad + \Bbar_N \Expectation{\Utilde_N \xtilde_0^{\TR}} \Abar_N^{\TR} + \Abar_N \Expectation{\xtilde_0\Utilde_N^{\TR}}\Bbar_N^{\TR}.
\end{align}

Applying the controller~\eqref{eqDiscreteLinearPolicy} results in the final covariance given by
\begin{equation}\label{eqFinalCovarianceInDiffusionlessCase}
\mSigma_{N+1} = \br{ \Abar_N + \Bbar_N \L}\covS\br{ \Abar_N + \Bbar_N \L}^{\TR} = \covF.
\end{equation}

The following proposition describes the diffusionless linear discrete covariance steering control algorithm:
\vspace{1em}

\begin{proposition}\label{lmDiffusionlessCovarianceSteering}
	 Let the controllable system~\eqref{eqTheSystemDiscrete}, with zero diffusion, and positive definite initial state covariance $\posdef{\covS}$, and let
	\begin{equation}\label{eqSVDsdefs}
	\V_0 \S_0 \V_0^{\TR} = \covS, \quad
	\V_F \S_F \V_F^{\TR} = \covF, \quad
	U_{\mOmega} \S_{\mOmega} \V_{\mOmega}^{\TR} = \mOmega,
	\end{equation}
be the singular value decompositions (SVDs) of the respective matrices, where
	\begin{equation}\label{eqOmegaDef}
	\mOmega \dfn {\S_F}^{\half}\V_F^{\TR}\invb{\Bbar_N\Bbar_N^{\TR}}  \Abar_N \V_0 {\S_0}^{\half}.
	\end{equation}

Then the optimal control gain $\L \in \Reals{(N r) \times n}$ that minimizes \eqref{eqCovarianceCostFunctionLinearController} subject to a constraint $\mSigma_{N+1} = \covF$, is given by
	\begin{equation}\label{eqDiscreteStateStatisticsCovNoDiffusionSolution}
		\Lstar = \Bbar_N^{\TR}\invb{\Bbar_N\Bbar_N^{\TR}}\br{\V_F {\S_F}^{\half}\U_{\mOmega} \V_{\mOmega}^{\TR}{\S_0}^{-\half}\V_0^{\TR} - \Abar_N}.
	\end{equation}
\end{proposition}
\vspace{1em}

\begin{proof}
Please see the appendix.
\end{proof}

\vspace{1em}


The proof of Proposition~\ref{lmDiffusionlessCovarianceSteering} reveals that the optimal control can also be obtained from 
	\begin{align}\label{eqSolutionForLUsingLambdaO}
	\Lstar &= - \Bbar_N^{\TR}\Lambda\invb{\I + \Bbar_N\Bbar_N^{\TR}\Lambda}\Abar_N,
	\end{align}
where  $\Lambda$ is the solution of a matrix Riccati equation.
\vspace{1em}

\begin{proposition}
	The matrix $\Lambda$ in~\eqref{eqSolutionForLUsingLambdaO} that satisfies the constraint~\eqref{eqFinalCovarianceInDiffusionlessCase}, and minimizes the cost function~\eqref{eqCovarianceCostFunctionLinearController}, satisfies the matrix Riccati equation
\begin{equation}\label{eqDiffusionelssRiccatiEquation}
\hspace*{-2mm}\br{\Theta\covF} \Lambda + \Lambda\br{\Theta\covF}^{\TR} + \Lambda\covF\Lambda  + \Theta\br{\covF - \Abar_N \covS \Abar_N^{\TR}}\Theta = 0,
\end{equation}
where $\Theta = \invb{\Bbar_N \Bbar_N^{\TR}}$.
\end{proposition}

\vspace{1em}

\begin{proof}

	Substituting $\L$ from~\eqref{eqSolutionForLUsingLambda} into the constraint~\eqref{eqRewrittenConstraintforCov}, and using matrix inversion identity, yields
	\begin{align}\label{eqConstExpanded}
	\covF &= \invb{\I + \Bbar_N \Bbar_N^{\TR}\Lambda}\Abar_N \covS \Abar_N^{\TR} \invb{\I + \Lambda \Bbar_N \Bbar_N^{\TR}},
	\end{align}
	which can be rewritten as~\eqref{eqDiffusionelssRiccatiEquation}.
\end{proof}

Note that the previous approach can be generalized to the case where it is required that the final covariance is only partially constrained,  \ie, given $\D \in \Reals{n_p \times n}$ with $n_p \leq n$ and final \textit{partial} covariance matrix $\covF \in \Reals{n_p \times n_p}$, the boundary condition for the state covariance at step $N+1$ is defined as
\begin{equation}\label{eqPartinalFinalCovariance}
\D \Expectation{\xtilde_{N+1} \xtilde_{N+1}^{\TR}} \D^{\TR} = \D \Sigma_{N+1} \D^{\TR} = \covF.
\end{equation}
Rewriting the above equation for a linear controller gain yields
\begin{equation}\label{eqDiscreteStateStatisticsCovNoDiffusionAndDemeanedQuadraticPartial}
\D \br{\Bbar_N \L + \Abar_N}\covS\br{\Bbar_N \L + \Abar_N}^{\TR} \D^{\TR}=\covF,
\end{equation}
which can be seen as the covariance-steering for diffusion-less system having transition matrices $\D\Bbar_N$ and $\D\Abar_N$, with the solution given by Proposition~\ref{lmDiffusionlessCovarianceSteering}.

\subsubsection{Non-zero Diffusion Case}\label{sec:subsubsec:GeneralCaseWithDiffusion}

Consider now the complete system given by~\eqref{eqTheSystemDiscrete}, including the diffusion term ($\G_k \neq 0$).

The

system~\eqref{eqTheSystemDiscrete} at time step $N+1$ can be viewed as a sum of $N+1$ uncorrelated ($\Expectation{\x^{(i)}_k {\x^{(j)}_m}^{\TR}} = 0,\ k, m, i, j \in [0, N+1],\ i\neq j$), diffusion-less sub-systems as follows
\begin{equation}\label{eqXksummation}
\x_{N+1} = \Sum{i=0}{N}{\x_{N+1}^{(i)}} + \G_N\w_N,
\end{equation}
where $\x_{N+1}^{(i)}$ for all $i=0,\ldots,N$ are computed, for all $k\in [i, N]$, from
\begin{equation}\label{eqSumOfSystems}
	\x^{(i)}_{k+1} = \A_k \x^{(i)}_{k} + \B_k \u^{(i)}_{k}, \quad \x^{(i)}_{i} = \left\{ \begin{matrix*}[l]
	\x_0, &\textrm{ for } i = 0, \\
	\G_{i-1}\w_{i-1}, &\textrm{ otherwise}, 
	\end{matrix*}\right.
\end{equation}
and $\x^{(i)}$ and $\u^{(i)}$ denote the state and the input of the $i$'th sub-system.
The final state can therefore be expressed as
\begin{gather}\label{eqStateAsACollection}
	\x_{N+1} = \Abar_N\x_0 + \Bbar_N\U^{(0)}_{0,N}  + \Sum{i=1}{N}{\Abar_{N, i}\G_{i-1}\w_{i-1} + \Bbar_{N, i}\U^{(i)}_{i, N}} + \G_N\w_N,
\end{gather}
where,
\begin{equation}\label{key}
\U^{(i)}_{k_1, k_2} \dfn \bmat{\u^{(i)}_{k_1} \\ \u^{(i)}_{k_1+1} \\ \vdots \\  \u^{(i)}_{k_2}}, \quad 0\leq k_1\leq k_2\leq N.
\end{equation}

We assume control laws with a linear dependence on $x^{(i)}_i$, that is, similarly to~\eqref{eqDiscreteLinearPolicy}, we
let $\L^{(k)} \in \Reals{(m(N-k+1)) \times n}$, $k\in[0,N]$, be a set of matrices, such that
\begin{equation}\label{eqLinearControlComponents}
\U^{(i)}_{i, N} = \begin{cases}
\L^{(i)} x^{(i)}_i, & i \in[1,N], \\
\L^{(0)} x_0 + \Expectation{\U_N}, & i = 0.
\end{cases}
\end{equation}
Since all states $x^{(i)}$ for $i \in[1,N]$ have zero mean, the mean of $x_{N+1}$ is governed by equation~\eqref{eqDiscreteStateStatisticsMean}. The covariance of the final state derived from~\eqref{eqStateAsACollection} is then given by
\begin{align}\label{eqFinalCovarianceWithDrift}
	& \Sigma_{N+1} = (\A_N + \Bbar_N\L^{(0)})\covS(\Abar_N + \Bbar_N\L^{(0)})^{\TR}  \notag\\
	&\quad + \Sum{i=1}{N}{(\A_{N, i} + \Bbar_{N, i}\L^{(i)})\G_{i-1}\G_{i-1}^{\TR} (\A_{N, i} + \Bbar_{N, i}\L^{(i)})^{\TR}} \notag\\
	&\quad + \G_N\G_N^{\TR}.
\end{align}

\begin{theorem}\label{thm:CovSteering}
	Let the system~\eqref{eqTheSystemDiscrete}, initial and final state means $\muS$ and $\muF$, and initial and final state covariance matrices $\possemidef{\covS}$ and $\possemidef{\covF}$.
Let $y_0 = x_0 - \muS$ and define, for $k\in[0,N]$,
	\begin{equation}\label{eqYkDef}
	 y_k = x_k - (\A_{k-1}\x_{k-1} + \B_{k-1}\u_{k-1}).
	\end{equation}
Furthermore, let $\Phi_k \in \Reals{n \times n }$ be given by
	\begin{equation}\label{eqPhiDefDrift}
	\Phi_k = \invb{\I + \Bbar_{N, k}\Bbar_{N, k}^{\TR}\Lambda}\A_{N, k},
	\end{equation}
	where $\Lambda = \Lambda^{\TR} \in \Reals{n \times n }$ is the solution of the matrix equation
	\begin{equation}\label{eqClosedLoopCovarianceEquivalence}
		\Sum{k=1}{N}{\Phi_k \G_{k-1}\G_{k-1}^{\TR}\Phi_k^{\TR}} +\Phi_0\covS\Phi_0^{\TR} = \covF - \G_{N}\G_{N}^{\TR}.
	\end{equation}

The optimal linear control law that minimizes the cost function~\eqref{eqCovarianceCostFunctionLinearController} subject to a constraints $\mSigma_{N+1} = \covF$ and $\mean_{N+1} = \muF$, and with the initial state mean $\muS$ and covariance $\covS$, is given by
	\begin{equation}\label{eqUkDriftCase}
	\u_k = \B_{N, k}^{\TR}\invb{\Bbar_N \Bbar_N^{\TR}}\br{\muF - \Abar_N\muS} + \Sum{i=0}{k}{\L^{(i)}_{k} y_i},
	\end{equation}
	where,
	\begin{equation}\label{eqControlGainForDriftCase}
	\L^{(i)}_{k} = -\B_{N, k}^{\TR}\Lambda\Phi_i.
	\end{equation}
\end{theorem}
\vspace{1em}

\begin{proof}
	Since the mean of the state is governed by~\eqref{eqDiscreteStateStatisticsMean}, the mean-steering solution $\Expectation{\U_N}$ is given by Proposition~\ref{lemma-steering-the-mean}, equation~\eqref{eqDiscreteOptimalSteeringMean}.
	
	The second part of the controller, $\Utilde_N$, having the linear form~\eqref{eqDiscreteLinearPolicy} is directed to minimizing the covariance-related cost~\eqref{eqCovarianceCostFunctionLinearController}, while adhering to the constraint $\Expectation{\xtilde_{N+1}\xtilde_{N+1}^{\TR}} = \covF$.
	
	The Lagrangian of the minimization problem~\eqref{eqCovarianceCostFunctionLinearController} subject to the constraint~\eqref{eqBoundaryConditions} is given by
	\begin{equation}\label{key}
	\Lagrangian(u, \Lambda) = \Trace{\Expectation{\Utilde_N \Utilde_N^\TR}} + \Trace{\Lambda(\Expectation{\xtilde_{N+1}\xtilde_{N+1}^{\TR}} - \covF)}.
	\end{equation}
	Using~\eqref{eqLinearControlComponents}, and~\eqref{eqFinalCovarianceWithDrift}, the Lagrangian can be rewritten in terms of $\L^{(i)}$, $i \in [0,N]$ as follows
	\begin{align}\label{key}
	&\Lagrangian(u, \Lambda) = \Tr\big\{\L^{(0)} \covS (\L^{(0)})^{\TR} - \Lambda\covF \notag\\
	&\qquad + \Lambda(\Abar_N + \Bbar_N\L^{(0)})\covS(\Abar_N + \Bbar_N\L^{(0)})^{\TR} \notag\\
	&\qquad + \sum_{i=1}^{N}\L^{(i)} \G_{i-1}\G_{i-1}^{\TR} (\L^{(i)})^{\TR} \notag\\
	&\qquad \quad	+ \Lambda(\A_{N, i} + \Bbar_{N, i}\L^{(i)})\G_{i-1}\G_{i-1}^{\TR} (\A_{N, i} + \Bbar_{N, i}\L^{(i)})^{\TR} \big\},
	\end{align}
	 yielding the following first and second order necessary conditions for a minimizer
	\begin{equation}\label{key}
	\L^{(i)} + \Bbar_{N,i}^{\TR}\Lambda\br{\A_{N,i} + \Bbar_{N,i} \L^{(i)}} = 0,
	\end{equation}
	and
	\begin{equation}\label{key}
	\posdef{\I + \Bbar_{N,i}^{\TR}\Lambda\Bbar_{N,i}}.
	\end{equation}
	Following a similar derivation as in Proposition~\ref{lmDiffusionlessCovarianceSteering}, the resulting optimal control gain is given by~\eqref{eqControlGainForDriftCase}, with $\Phi_k$ given by~\eqref{eqPhiDefDrift}. Substituting this control back into the constraint equation~\eqref{eqFinalCovarianceWithDrift}, results in the closed-loop covariance equation~\eqref{eqClosedLoopCovarianceEquivalence}.
	
	Therefore, the matrix $\Lambda$ that satisfies the constraint~\eqref{eqClosedLoopCovarianceEquivalence} provides the optimal gains for the optimal controller~\eqref{eqUkDriftCase}.
\end{proof}
\vspace{1em}

Note that the controller in~\eqref{eqUkDriftCase} can be efficiently calculated by updating the vector $\U_{k,N}$ at every step $k$ (starting from $k=0$) by
\begin{align}\label{eqEfficient1}
	\U_{0,N} &=  \Bbar_N^{\TR}\invb{\Bbar_N \Bbar_N^{\TR}}\br{\muF - \Abar_N\muS}, \notag\\
	\U_{k,N} &=  \U_{k,N} + \L^{(k)} \br{x_k - \hat{x}_k},
\end{align}
where
\begin{align}\label{eqEfficient2}
	\hat{x}_0 = \muS, \qquad \hat{x}_{i+1} = \A_i \hat{x}_i + \B_i u_i.
\end{align}

The non-negativity of the left-hand side of~\eqref{eqClosedLoopCovarianceEquivalence} yields
\begin{equation}\label{key}
\possemidef{\covF - \G_{N}\G_{N}^{\TR}},
\end{equation}
which is exactly the condition for solvability for the covariance steering problem provided in~\cite[Proposition 5.1]{Beghi1996}.

\section{Relation with LQG}

The stochastic control problem formulated in Section~\ref{sec:secProblemStatement:ProblemFormulation} can be also viewed as a special case of the standard discrete {LQG}~\cite[p.264]{Astrom1970}. This similarity will be detailed in this section, focusing on the covariance control, thus assuming a zero-mean state.

A LQG control problem is formulated as follows\footnote{Here a shortened version of the LQG problem is given. The original version~\cite[p.264]{Astrom1970} includes quadratic cost matrices for the control effort, and for the state at each step.}. 
Given a stochastic discrete linear system~\eqref{eqTheSystemDiscrete}, with the noise described by~\eqref{eqNoiseMdl1} and~\eqref{eqNoiseMdl2}, and the initial state drawn from a normal distribution
$\x_0 \sim \normal{0,\covS}$,
let $\Q_f \in \Reals{n \times n}$ be a symmetric matrix. 
The optimal control that minimizes the cost function
\begin{equation}\label{eqLQRcostFunction}
J(\u_0, ..., \u_N) = \Expectation{\Sum{k=0}{N}{ \u_{k}^\TR  \u_{k}} + \x_{N+1}^{\TR} \Q_f \x_{N+1}},
\end{equation}
subject to the~\eqref{eqTheSystemDiscrete} is given by
\begin{equation}\label{eqLQRcontroller1}
\u_k = -\L_k \x_k,
\end{equation}
where
\begin{equation}\label{eqLQRcontroller2}
\L_k = \invb{\I + \B_k^{\TR} \P_{k+1}\B_k} \B_k^{\TR} \P_{k+1} \A_k,
\end{equation}
and $\P_k=\P_k^{\TR}$ is given by the backward-recursive equation
\vspace{1em}
\begin{theorem}
	Let the system~\eqref{eqTheSystemDiscrete}, with zero-mean states, and initial and final state covariance matrices $\covS$ and $\covF$. Let $\Q_f \in \Reals{n \times n}$ be a symmetric matrix. Assume that the LQG controller that minimizes the cost function
	\begin{equation}\tag{\ref{eqLQRcostFunction}}
	J(\u_0, ..., \u_N) = \Expectation{\Sum{k=0}{N}{ \u_{k}^\TR  \u_{k}} + \x_{N+1}^{\TR} \Q_f \x_{N+1}},
	\end{equation}
	subject to the dynamics~\eqref{eqTheSystemDiscrete}, results in the final state covariance being equal to $\covF$. 
Then, this controller coincides with the optimal controller given by the problem described in  Theorem~\ref{thm:CovSteering}, that is, it minimizes \eqref{eqDiscreteCostFunction} subject to the dynamics~\eqref{eqTheSystemDiscrete} and the boundary constraints
	\begin{equation}\label{eqBndCondsInProof}
	\x_0 \sim \normal{0,\covS}, \qquad \x_{N+1} \sim \normal{0,\covF},
	\end{equation}
	with $\Lambda = \Q_f$.
\end{theorem}
\vspace{1em}

\begin{proof}
The Lagrangian can be written as
\begin{align}\label{eqLagrangianInProofOfLQRtrm}
\Lagrangian &= \Expectation{\Sum{k=0}{N}{ \u_{k}^\TR  \u_{k}}} + \Trace{\Lambda\br{\Sigma_{N+1} - \covF}} \notag\\
	&= \Expectation{\Sum{k=0}{N}{ \u_{k}^\TR  \u_{k}} + \x_{N+1}^{\TR}\Lambda\x_{N+1}} - \Trace{\Lambda\covF}.
\end{align}
Given that $\Lambda = \Q_f$, minimizing the Lagrangian~\eqref{eqLagrangianInProofOfLQRtrm} yields the same result as minimizing
\begin{equation}\label{key}
\Lagrangian = \Expectation{\Sum{k=0}{N}{ \u_{k}^\TR  \u_{k}} + \x_{N+1}^{\TR}\Q_f\x_{N+1}},
\end{equation}
and the optimal solution is given by the LQG controller described by equations~\eqref{eqLQRcontroller1}-\eqref{eqPequationLQR}. Since, by construction, this solution agrees with the boundary conditions~\eqref{eqBndCondsInProof}, it is also a solution of the covariance steering problem.
\end{proof}
\vspace{1em}
\begin{corollary}
	Assume $\Lambda$, which solves the optimal control problem given by~\eqref{eqDiscreteCostFunction},
	is unique.
	Then, the controller~\eqref{eqUkDriftCase} coincides with the LQG controller that minimizes the cost function:
	\begin{equation}\label{key}
	J(\u_0, ..., \u_N) = \Expectation{\Sum{k=0}{N}{ \u_{k}^\TR  \u_{k}} + \x_{N+1}^{\TR} \Lambda \x_{N+1}}
	\end{equation}
\end{corollary}
\vspace{1em}

\begin{proof}
	Recall that the Lagrangian of the optimal control problem given by~\eqref{eqDiscreteCostFunction} can be written as~\eqref{eqLagrangianInProofOfLQRtrm}, namely
	\begin{equation*}
	\Lagrangian = \Expectation{\Sum{k=0}{N}{ \u_{k}^\TR  \u_{k}} + \x_{N+1}^{\TR}\Lambda\x_{N+1}} - \Trace{\Lambda\covF}.
	\end{equation*}
	Since $\Lambda=\Lambda^{\TR}$ is given,
	\begin{align}\label{eqLagrang2}
	U_N 
		& = \arg \min_{\U_N} \Expectation{\Sum{k=0}{N}{ \u_{k}^\TR  \u_{k}} + \x_{N+1}^{\TR}\Lambda\x_{N+1}} - \Trace{\Lambda\covF} \notag\\
		&= \arg \min_{\U_N} \Expectation{\Sum{k=0}{N}{ \u_{k}^\TR  \u_{k}} + \x_{N+1}^{\TR}\Lambda\x_{N+1}},
	\end{align}
	subject to the dynamics~\eqref{eqTheSystemDiscrete}. The solution to~\eqref{eqLagrang2} is given by the LQG controller~\eqref{eqLQRcontroller1}, and minimizes the cost~\eqref{eqLQRcostFunction} with $\Q_F = \Lambda$.
\end{proof}
\vspace{1em}

The presented similarity to the LQG case allows an iterative implementation of the covariance steering controller, using equations~\eqref{eqLQRcontroller1}-\eqref{eqPequationLQR}.

Note that the presented results coincide with the results in~\cite{Beghi1996}. In fact, equation~(5.5) in~\cite{Beghi1996} is exactly equation~\eqref{eqClosedLoopCovarianceEquivalence}, with the right closed-loop transition matrices.

\section{Numerical Example}\label{secNumericalExample}

In this section the performance of the algorithm is tested using two simple examples: First, a second order LTI system, and then a fourth-order linear time-varying system, which is derived from linearizing and discretizing a non-linear cart-pole dynamics along a particular trajectory are used to demonstrate the results of the previous sections.

\subsection{Example: Linear Time Invariant System}

Consider the LTI stochastic discrete system:
\begin{gather}\label{eqLinearSystemExample}
\A = \bmat{1.9986 & -1 \\ 1 & 0}, \quad \B = \bmat{0.03125 \\ 0}, \quad \G = \bmat{0 \\ 0.03},
\end{gather}
with $N=100$, and the boundary conditions
\begin{equation}
\begin{aligned}
&\muS = \bmat{-1 \\ 1}, & \qquad &\covS = \bmat{2 & -1\\-1 & 3}, \\
&\muF = \bmat{1 \\ -1}, & \qquad &\covF = \bmat{1 & 0.1\\0.1 & 2}.
\end{aligned}
\end{equation}

Solving for $\Lambda$ in equation~\eqref{eqClosedLoopCovarianceEquivalence} yields:
\begin{equation}\label{key}
\Lambda = \bmat{12.225
	& -3.398
	\\ -3.398
	     &     6.543
	     }.
\end{equation}
This value was used in the controller form~\eqref{eqUkDriftCase}, as well as in the LQG form, by setting $\P_{N+1} = \Lambda$ in~\eqref{eqPequationLQR}. In addition, the mean-steering open-loop controller was calculated using equation~\eqref{eqDiscreteOptimalSteeringMean}.

\begin{figure}[thpb]
	\centering
	\includegraphics[width=0.95\linewidth]{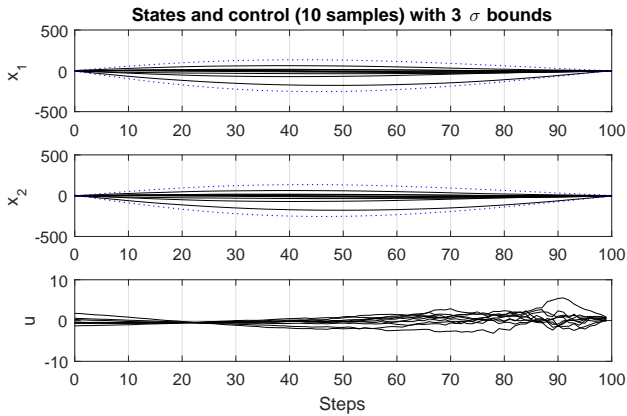}
	\caption{LTI model: states and controls from 10 randomly-chosen trajectories - Covariance Controller algorithm}
	\label{LTI_Mdl_StatesAndControlSample}
\end{figure}

\begin{figure}[thpb]
	\centering
	\includegraphics[width=0.95\linewidth]{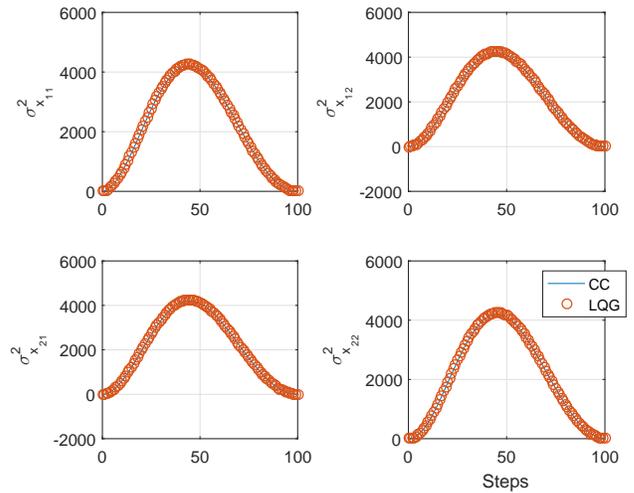}
\caption{LTI model: State covariance, Covariance Controller (CC) vs. LQG}
\label{LTI_Mdl_StateCovariance}
\end{figure}

\begin{figure}[thpb]
	\centering
	\includegraphics[width=0.95\linewidth]{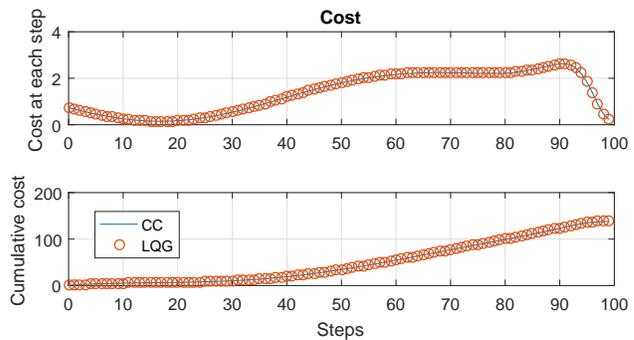}
	\caption{LTI model: Control costs, Covariance Controller (CC) vs. LQG}
\label{LTI_Mdl_Cost}
\end{figure}

As expected, the two algorithms give exactly the same result, in terms of the optimal trajectories and control gains. The results are depicted on Figures~\ref{LTI_Mdl_StatesAndControlSample}--\ref{LTI_Mdl_Cost}. The statistical values on the graphs were generated using Monte-Carlo simulations (20,000 runs) for the Covariance Control algorithm, and by using analytical expressions for the LQG algorithm.


\subsection{Example: Linear Time-Varying System}

In this example we consider a time-varying linear system generated by the linearization of a nonlinear cart-pole model around a nominal trajectory.
Let $y$ denote the cart's position, let $u$ denote the force pushing the cart, and let $\theta$ denote the pole's angle measured from vertical axis so that $\theta=0$ indicates the configuration when
the pole points vertically downwards. 
The equations of the of the cart-pole system are
\begin{align}\label{eqCartPoleEquationsOfMotion}
	\ddot{\theta} &= \frac{-\br{u+m_p l \dot{\theta}^2\sin \theta} \cos \theta  - \br{m_c+m_p}g\sin \theta}{l\br{m_c+m_p \sin^2 \theta}}, \notag\\
	\ddot{y} &= \frac{u+m_p \sin \theta \br{l \dot{\theta}^2+g \cos \theta}}{m_c+m_p \sin^2 \theta}.
\end{align}
The following parameters were used in the numerical simulations: $m_p = 0.01 [\textrm{kg}], m_c = 1 [\textrm{kg}], \ell = 0.25 [\textrm{m}], g = 9.81 [{\textrm{m}}/{\textrm{sec}^2}]$.
The equations of motion were linearized about a trajectory that brings the pole rom the downward position $\theta_0 = 0$ to the upward position $\theta_F = \pi$ in 1 second 
, and then discretized using Euler's method with sampling time of $T_s=0.001~\mathrm{sec}$, resulting in a linear discrete time-varying system with states defined as
\begin{equation}\label{key}
\x \dfn \bmat{\delta \theta & \delta \dot{\theta} & \delta y & \delta \dot{y}}^{\TR},
\end{equation}
where $\delta \theta$, $\delta \dot{\theta}$, $\delta y$, and $\delta \dot{y}$ denote deviations from the nominal values of $\theta$, $\dot{\theta}$, $y$, and $\dot{y}$ respectively. To this model, a disturbance noise was added, with:
\begin{equation}\label{key}
\G = \bmat{ 0           &           0.004            &             0             &         0.008	}^{\TR}.
\end{equation}

As in the previous example, it is assumed that the state is fully measured.
The initial and the final states are chosen as
\begin{equation}\label{eqNumExampleCartPoleBCs}
\begin{aligned}
	\muS = \muF = 0_{4 \times 1},\quad \covS = \covF = \diag { 0.01,0.01,0.01,0.01},\\
\end{aligned}
\end{equation}

\begin{figure}[thbp]
	\centering
	\includegraphics[width=0.95\linewidth]{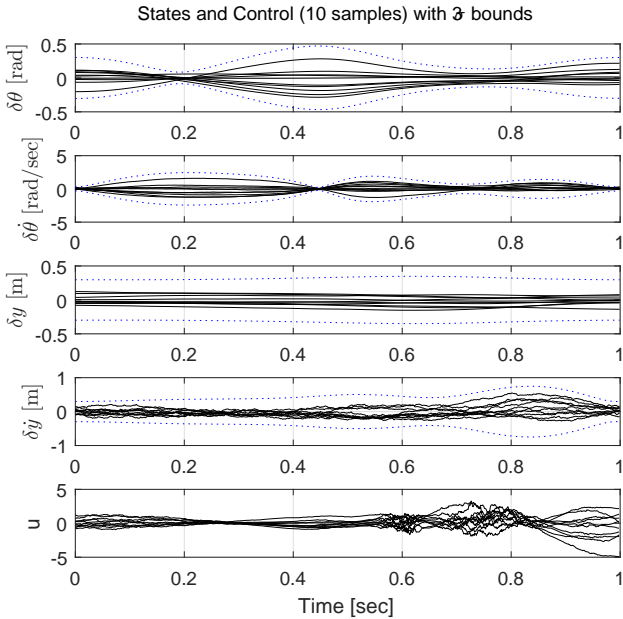}
	\caption{Linearized Cart-pole: 10 samples of the states and the controls.}
	\label{CartPole_StatesAndControlSample}
\end{figure}

\begin{figure}[thbp]
	\centering
	\includegraphics[width=0.95\linewidth]{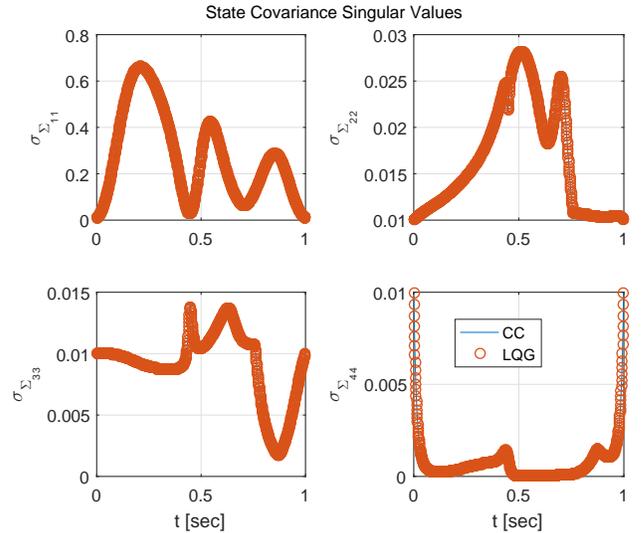}
	\caption{Linearized Cart-pole: State covariance Singular Values.}
	\label{CartPole_StateCovarianceSVs}
\end{figure}

\begin{figure}[thbp]
	\centering
	\includegraphics[width=0.95\linewidth]{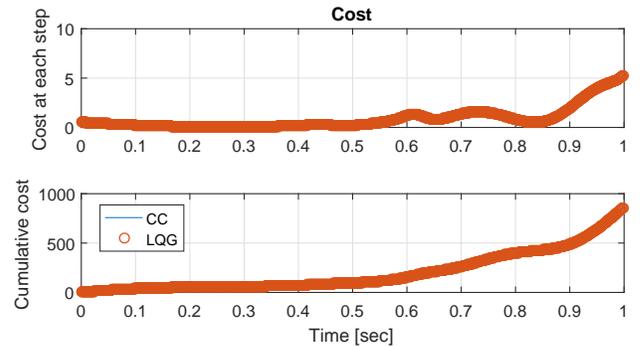}
	\caption{Linearized Cart-pole: Control cost.}
	\label{CartPole_Cost}
\end{figure}

The results are shown in Figures~\ref{CartPole_StatesAndControlSample}-\ref{CartPole_Cost}. Figure~\ref{CartPole_StatesAndControlSample} exhibits 10 randomly-generated closed-loop trajectories (states and control), and the $3\sigma$ bounds calculated from 20,000 Monte-Carlo runs. The controller costs are shown in Figure~\ref{CartPole_Cost}. Figure~\ref{CartPole_StateCovarianceSVs} depicts evaluation of state covariance singular values through time. 


Similarly to the LTI example, the results show that the two algorithms give exactly the same results.

\section{Conclusions}

In this work we have derived a minimum-control-effort optimal steering solution for linear time-varying discrete stochastic systems, subject to boundary conditions in the form of Gaussian distribution parameters. Having presented the influence of the diffusion at each time-step on the final covariance, we have formulated a condition for calculating the optimal control law from a class of linear-state-dependent control laws. 
The resulting controller set consists of ``open-loop'' inputs, which are recalculated at each step based on the diffusion term reconstruction from the previous step.

In addition, we have shown that the solution to the covariance steering problem coincides with the solution to a specially-formulated LQG problem. This similarity allowed an efficient calculation of the controller values using a backward-propagated discrete-time dynamic Riccati equation, as well as a justification for using a linear feedback controller for the covariance steering.

We have demonstrated the algorithm performance using simple numeral examples, showing that the covariance steering algorithm coincides with the respective LQR solution.

Future work will address the conditions for the existence of the solution, and the algorithm applicability for the covariance steering of non-linear systems.

\section*{Acknowledgment}

This work was supported by NSF award NRI-1426945 and ARO award W911NF-16-1-0390.

\bibliographystyle{IEEEtran}
\bibliography{Covariance-control-methods}


\appendix
\section*{Proof of Proposition~\ref{lmDiffusionlessCovarianceSteering}}
	Let $\Xi \dfn \Abar_N + \Bbar_N \L$. Then the constraint~\eqref{eqFinalCovarianceInDiffusionlessCase} can be written as
	\begin{equation}\label{eqRewrittenConstraintforCov}
	\Xi \covS \Xi^{\TR} = \covF.
	\end{equation}
	First we show feasibility.
	Substituting~\eqref{eqDiscreteStateStatisticsCovNoDiffusionSolution} into \eqref{eqRewrittenConstraintforCov} yields
	\begin{align}\label{eqXiEquationSubstituted}
	\Xi &= \Abar_N + \Bbar_N\Bbar_N^{\TR}\invb{\Bbar_N\Bbar_N^{\TR}}\br{\V_F {\S_F}^{\half}\U_{\mOmega} \V_{\mOmega}^{\TR}{\S_0}^{-\half}\V_0^{\TR} - \Abar_N} \notag\\
	&= V_F {\S_F}^{\half}\U_{\mOmega} \V_{\mOmega}^{\TR}{\S_0}^{-\half}\V_0^{\TR},
	\end{align}
	and hence
	\begin{align}\label{key}
	\Xi \covS \Xi^{\TR} &= V_F {\S_F}^{\half}\U_{\mOmega} \V_{\mOmega}^{\TR} {\S_0}^{-\half} \V_0^{\TR} \V_0 \S_0 \V_0^{\TR} \V_0 {\S_0}^{-\half} \V_{\mOmega} \U_{\mOmega}^{\TR} {\S_F}^{\half} V_F^{\TR} \notag\\
	&=V_F {\S_F} V_F^{\TR} = \covF.
	\end{align}
	
	To show optimality, introduce the Lagrangian of the equality constraint minimization problem~\eqref{eqCovarianceCostFunctionLinearController} and~\eqref{eqFinalCovarianceInDiffusionlessCase}
	\begin{align}
	\Lagrangian(\L, \Lambda) &= \Trace{\L \covS \L^{\TR}} + \Trace{\Lambda \big( \Xi\covS\Xi^{\TR} - \covF \big) }
	\end{align}
	where $\Lambda \in \Reals{n \times n}$.
	Without loss of generality we assume that $\Lambda = \Lambda^{\TR}$. The first-order optimality condition $\Lagrangian_{\L}(\L, \Lambda) = 0$ yields:
	\begin{equation}\label{eqLargarngianFirstOrder}
	\L + \Bbar_N^{\TR}\Lambda\br{\Abar_N + \Bbar_N \L} = \L + \Bbar_N^{\TR}\Lambda\Xi= 0,
	\end{equation}
	whereas the second order condition $\Lagrangian_{\L \L}(\L, \Lambda) = 0$ yields
	\begin{equation}\label{eqLargarngianSecondOrder}
	\posdef{\I + \Bbar_N^{\TR}\Lambda\Bbar_N}.
	\end{equation}
	It follows that
	\begin{align}\label{eqSolutionForLUsingLambda}
	\L &= - \Bbar_N^{\TR}\Lambda\invb{\I + \Bbar_N\Bbar_N^{\TR}\Lambda}\Abar_N.
	\end{align}
	Substituting this value of $L$ into the constraint \eqref{eqFinalCovarianceInDiffusionlessCase} yields
	\begin{align}\label{eqXiExpanded}
	\Xi &= \Abar_N - \Bbar_N \Bbar_N^{\TR}\Lambda\invb{\I + \Bbar_N\Bbar_N^{\TR}\Lambda}\Abar_N \notag\\
	&= \invb{\I + \Bbar_N \Bbar_N^{\TR}\Lambda}\Abar_N.
	\end{align}
	Using the SVDs~\eqref{eqSVDsdefs} we can rewrite the constraint~\eqref{eqRewrittenConstraintforCov} as 
	\begin{equation}\label{eqExpressionWithR}
	\Xi \V_0 \S_0^{\half}\R^{\TR} = \V_F\S_F^{\half},
	\end{equation}	
	where $\R$ is an orthogonal matrix. 
	%
	Combining~\eqref{eqExpressionWithR} with~\eqref{eqXiExpanded} yields
	\begin{equation}\label{key}
	\Bbar_N \Bbar_N^{\TR}\Lambda = \Abar_N\V_0^{\TR}\S_0^{\half}\R^{\TR}\S_F^{-\half}\V_F^{\TR} - \I,
	\end{equation}
	and the resulting optimal gain is
	\begin{equation}\label{key}
	\L^{\star} = \Bbar_N^{\TR}\invb{\Bbar_N\Bbar_N^{\TR}}\br{\V_F {\S_F}^{\half}\R{\S_0}^{-\half}\V_0^{\TR} - \Abar_N}.
	\end{equation}
	
	In order to find $\R$, the optimal gain equation is substituted into the cost function $J_{\Sigma}$, resulting in
	\begin{align}\label{eqSigmaCostSubstituted}
	J_{\Sigma} &= \Trace{\Bbar_N^{\TR}\invb{\Bbar_N \Bbar_N^{\TR}}\br{\V_F\S_F^{\half}\R\S_0^{-\half}\V_0^{\TR} - \Abar_N} \V_0 \S_0 \V_0^{\TR}{\L^{\star}}^{\TR}} \notag\\
	&= \Trace{\invb{\Bbar_N\Bbar_N^{\TR}} \br{\covF + \Abar_N\covS\Abar_N^{\TR}}} - 2\Trace{\R^{\TR}\U_{\mOmega} \S_{\mOmega} \V_{\mOmega}^{\TR}}
	\end{align}
	where $\mOmega$ was defined in~\eqref{eqOmegaDef}. The minimum of the cost~\eqref{eqSigmaCostSubstituted} is attained by maximizing the term $\Trace{\R^{\TR}\U_{\mOmega}}$, yielding
	\begin{align}
	\R^{\star} = \arg\min_{\R\in \mathcal{U}^n} J_{\Sigma}
	= \arg\max_{\R\in \mathcal{U}^n} \Trace{\R^{\TR}\U_{\mOmega} \S_{\mOmega} \V_{\mOmega}^{\TR}}
	= \U_{\mOmega} \V_{\mOmega}^{\TR},
	\end{align}
	where the last equation follows from the von Newmann trace inequality~\cite{Mirsky1975}.
	Substituting $\R^{\star}$ into the optimal gain $\L^{\star}$ yields~\eqref{eqDiscreteStateStatisticsCovNoDiffusionSolution}.


\end{document}

%% file: LatexCommon/command.tex
\usepackage{ifthen}
\usepackage{mathtools} 
\usepackage{relsize} 

\makeatletter
\newcommand*\rel@kern[1]{\kern#1\dimexpr\macc@kerna}
\newcommand*\widebar[1]{%
	\begingroup
	\def\mathaccent##1##2{%
		\rel@kern{0.8}%
		\overline{\rel@kern{-0.8}\macc@nucleus\rel@kern{0.2}}%
		\rel@kern{-0.2}%
	}%
	\macc@depth\@ne
	\let\math@bgroup\@empty \let\math@egroup\macc@set@skewchar
	\mathsurround\z@ \frozen@everymath{\mathgroup\macc@group\relax}%
	\macc@set@skewchar\relax
	\let\mathaccentV\macc@nested@a
	\macc@nested@a\relax111{#1}%
	\endgroup
}
\makeatother

\newcommand{\ie}{{i.e.}\xspace}

\newcommand{\dfn}{\triangleq}

\let\oldforall\forall
\renewcommand{\forall}{\ensuremath{\; \oldforall \;}}

\newcommand{\Sum}[3]{\sum\limits_{#1}^{#2} {#3} \hspace{2pt}}

\newcommand{\TR}{{\hspace{-1pt}{\top}\hspace{-1pt}}}

\newcommand{\diag}[1]{\mathrm{diag}\left[#1\right]}

\newcommand{\normal}[1]{\mathcal{N}(#1)}

\renewcommand{\vec}[1]{\ensuremath{{#1}}}
\newcommand{\mtx}[1]{\ensuremath{{#1}}}

\newcommand{\br}[1]{{\hspace{-1pt}(#1)}}

\newcommand{\brs}[1]{{[{#1}]}}

\newcommand{\brf}[1]{{\left\{{#1}\right\}}}

\newcommand{\Expectation}[1]{\mathbb{E}[#1]}
\newcommand{\Rfield}{\mathbb{R}}
\newcommand{\field}[1]{\Rfield^{#1}}

\newcommand{\Reals}[1]{\field{#1}}

\newcommand{\bmat}[1]{\begin{bmatrix} #1 \end{bmatrix}}

\newcommand{\invb}[1]{\br{#1}^{\hspace{-2pt}-1}}

\newcommand{\I}[0]{\mathsf{I}}

\newcommand{\half}{\frac{1}{2}}

\DeclareMathOperator{\Tr}{Tr}

\newcommand{\Trace}[1]{\Tr\brs{#1}} 

\newcommand{\posdef}[1]{{#1} \succ 0}
\newcommand{\possemidef}[1]{{#1} \succeq 0}

%% file: MacrosAndDefinitions.tex
\newcommand{\atTime}[1]{\ifthenelse{\equal{#1}{}}{}{\br{#1}}}
\newcommand{\atStep}[1]{\ifthenelse{\equal{#1}{}}{}{_{#1}}}

\newcommand{\x}{\vec{x}}
\renewcommand{\u}{\vec{u}}
\newcommand{\w}{\vec{w}}

\newcommand{\W}{\vec{W}}

\newcommand{\xtilde}{\widetilde{\x}}
\newcommand{\Utilde}{\widetilde{\U}}

\renewcommand{\S}{\mtx{S}}
\newcommand{\V}{\mtx{V}}
\newcommand{\U}{\mtx{U}}
\newcommand{\R}{\mtx{R}}

\newcommand{\D}{\mtx{D}}
\renewcommand{\L}{\mtx{L}}

\newcommand{\Q}{\mtx{Q}}
\newcommand{\G}{\mtx{G}}

\newcommand{\A}{\mtx{A}}
\newcommand{\B}{\mtx{B}}

\newcommand{\Abar}{\widebar{\A}}
\newcommand{\Bbar}{\widebar{\B}}
\newcommand{\Gbar}{\widebar{\G}}

\renewcommand{\P}{\mtx{P}}

\newcommand{\mOmega}{\mtx{\Omega}}

\newcommand{\mSigma}{\mtx{\Sigma}}

\newcommand{\mean}{\vec{\mu}}
\newcommand{\muS}{\mean_0}
\newcommand{\covS}{\mSigma_0}

\newcommand{\muF}{\mean_F}
\newcommand{\covF}{\mSigma_F}

\newcommand{\optU}{\U^{\star}}

\newcommand{\Lstar}{{\L^{\star}}}

\newcommand{\Lagrangian}{\mathfrak{L}}

\newcommand{\deltaX}[1]{\delta\vx}
\newcommand{\deltaU}[1]{\delta\vu}

\newcommand{\scrFfun}[1]{\mathfrak{F}\atTime{k}}
\newcommand{\scrZfun}[1]{\mathfrak{Z}\atTime{k}}
\newcommand{\scrLfun}[1]{\mathfrak{L}\atTime{k}}

\renewcommand{\k}{\vec{k}}

















%% file: ms.bbl
\begin{thebibliography}{10}
\providecommand{\url}[1]{#1}
\csname url@samestyle\endcsname
\providecommand{\newblock}{\relax}
\providecommand{\bibinfo}[2]{#2}
\providecommand{\BIBentrySTDinterwordspacing}{\spaceskip=0pt\relax}
\providecommand{\BIBentryALTinterwordstretchfactor}{4}
\providecommand{\BIBentryALTinterwordspacing}{\spaceskip=\fontdimen2\font plus
\BIBentryALTinterwordstretchfactor\fontdimen3\font minus
  \fontdimen4\font\relax}
\providecommand{\BIBforeignlanguage}[2]{{%
\expandafter\ifx\csname l@#1\endcsname\relax
\typeout{** WARNING: IEEEtran.bst: No hyphenation pattern has been}%
\typeout{** loaded for the language `#1'. Using the pattern for}%
\typeout{** the default language instead.}%
\else
\language=\csname l@#1\endcsname
\fi
#2}}
\providecommand{\BIBdecl}{\relax}
\BIBdecl

\bibitem{Astrom1970}
K.~J. {\AA}str{\"o}m, \emph{Introduction to Stochastic Control Theory},
  R.~Bellman, Ed.\hskip 1em plus 0.5em minus 0.4em\relax Academic Press New
  York and London, 1970.

\bibitem{Shahrokhi2016}
S.~Shahrokhi, A.~Mahadev, and A.~T. Becker, ``Algorithms for shaping a particle
  swarm with a shared control input using boundary interaction,'' \emph{arXiv
  preprint arXiv:1609.01830}, 2016.

\bibitem{Vinante2008}
A.~Vinante, M.~Bignotto, M.~Bonaldi, M.~Cerdonio, L.~Conti, P.~Falferi,
  N.~Liguori, S.~Longo, R.~Mezzena, A.~Ortolan \emph{et~al.}, ``Feedback
  cooling of the normal modes of a massive electromechanical system to
  submillikelvin temperature,'' \emph{Physical Review Letters}, vol. 101,
  no.~3, p. 033601, 2008.

\bibitem{Hotz1987}
A.~Hotz and R.~E. Skelton, ``Covariance control theory,'' \emph{International
  Journal of Control}, vol.~46, no.~1, pp. 13--32, 1987.

\bibitem{GrigoriadisSkelton1997}
K.~M. Grigoriadis and R.~E. Skelton, ``Minimum-energy covariance controllers,''
  \emph{Automatica}, vol.~33, no.~4, pp. 569--578, 1997.

\bibitem{Xu1992}
J.-H. Xu and R.~Skelton, ``An improved covariance assignment theory for
  discrete systems,'' \emph{IEEE Transactions on Automatic Control}, vol.~37,
  no.~10, pp. 1588--1591, 1992.

\bibitem{Collins1987}
E.~Collins and R.~Skelton, ``A theory of state covariance assignment for
  discrete systems,'' \emph{IEEE Transactions on Automatic Control}, vol.~32,
  no.~1, pp. 35--41, 1987.

\bibitem{Chen.etal2016}
Y.~Chen, T.~T. Georgiou, and M.~Pavon, ``Optimal steering of a linear
  stochastic system to a final probability distribution, {Part} {I},''
  \emph{IEEE Transactions on Automatic Control}, vol.~61, no.~5, pp.
  1158--1169, May 2016.

\bibitem{Chen.etal2016a}
------, ``Optimal steering of a linear stochastic system to a final probability
  distribution, {Part} {II},'' \emph{IEEE Transactions on Automatic Control},
  vol.~61, no.~5, pp. 1170--1180, May 2016.

\bibitem{Chen2016Thesis}
Y.~Chen, ``Modeling and control of collective dynamics: From {Schr{\"o}dinger}
  bridges to optimal mass transport,'' Ph.D. dissertation, University of
  Minnesota, 2016.

\bibitem{Schrodinger1931}
E.~Schr{\"o}dinger, \emph{{\"U}ber die {Umkehrung} der {Naturgesetze}}.\hskip
  1em plus 0.5em minus 0.4em\relax Verlag Akademie der Wissenschaften in
  Kommission bei {Walter} de {Gruyter} u. {Company}, 1931.

\bibitem{Kantorovich1942}
L.~V. Kantorovich, ``On the transfer of masses,'' in \emph{Dokl. Akad. Nauk.
  SSSR}, vol.~37, no. 7-8, 1942, pp. 227--229.

\bibitem{Halder2016ACC}
A.~Halder and E.~D. Wendel, ``Finite horizon linear quadratic {Gaussian}
  density regulator with {Wasserstein} terminal cost,'' in \emph{American
  Control Conference}, Boston, MA, USA, July 6--8 2016, pp. 7249--7254.

\bibitem{Bakolas2016ACC}
E.~Bakolas, ``Optimal covariance control for stochastic linear systems subject
  to integral quadratic state constraints,'' in \emph{American Control
  Conference}, Boston, MA, USA, July 6--8 2016, pp. 7231--7236.

\bibitem{Bakolas2016CDC}
------, ``Optimal covariance control for discrete-time stochastic linear
  systems subject to constraints,'' in \emph{55th Conference on Decision and
  Control}.\hskip 1em plus 0.5em minus 0.4em\relax Las Vegas, NV, USA: IEEE,
  December 12-14 2016, pp. 1153--1158.

\bibitem{Beghi1996}
A.~Beghi, ``On the relative entropy of discrete-time {Markov} processes with
  given end-point densities,'' \emph{{IEEE} Transactions on Information
  Theory}, vol.~42, no.~5, pp. 1529--1535, 1996.

\bibitem{Hamedani2001}
G.~Hamedani and H.~Volkmer, ``Certain characterizations of normal distribution
  via transformations,'' \emph{Journal of Multivariate Analysis}, vol.~77,
  no.~2, pp. 286--294, 2001.

\bibitem{Brockett2007}
R.~Brockett, ``Optimal control of the {Liouville} equation,'' \emph{AMS IP
  Studies in Advanced Mathematics}, vol.~39, p.~23, 2007.

\bibitem{Corless2003}
M.~J. Corless and A.~Frazho, \emph{Linear Systems and Control: an Operator
  Perspective}.\hskip 1em plus 0.5em minus 0.4em\relax CRC Press, 2003.

\bibitem{Mirsky1975}
L.~Mirsky, ``A trace inequality of {John} {von} {Neumann},'' \emph{Monatshefte
  f{\"u}r Mathematik}, vol.~79, no.~4, pp. 303--306, 1975.

\end{thebibliography}
